\documentclass[11pt,a4paper]{article}
\usepackage[top=2.5cm,bottom=2.5cm,left=2.2cm,right=2.2cm]{geometry}
\usepackage[T1]{fontenc}
\usepackage[utf8]{inputenc}
\usepackage{color}
\usepackage{graphicx}
\usepackage{amsfonts}
\usepackage{extarrows}
\usepackage{amsmath,amsthm,amssymb,color}
\usepackage{hyperref}
\usepackage{eepic}
\usepackage{lineno}
\usepackage{enumerate}	
\usepackage{paralist}
\usepackage{cite}
\usepackage{algorithm}
\usepackage{algorithmicx}
\usepackage{algpseudocode}
\usepackage{authblk}
\usepackage{mathtools}
\usepackage{pifont}
\usepackage[numbers,sort&compress]{natbib}

\usepackage{tikz}

\newtheorem{theorem}{Theorem}[section]

\newtheorem{conjecture}[theorem]{Conjecture}
\newtheorem{problem}[theorem]{Problem}

\theoremstyle{definition}

\newtheorem{claim}{\indent Claim}

\newenvironment{wst}
{\setlength{\leftmargini}{1.5\parindent}
 \begin{itemize}
 \setlength{\itemsep}{-1.1mm}}
{\end{itemize}}

\begin{document}
	\title{\bf Constructions of  minimally $t$-tough regular graphs}
	\author{Kun Cheng\footnote{Email: chengkunmath@163.com.}}
	\author{Chengli Li\footnote{Email: lichengli0130@126.com.}}
	\author{Feng Liu\footnote{Email: liufeng0609@126.com.}}
	
	\affil{\footnotesize Department of Mathematics,
		East China Normal University, Shanghai, 200241, China}
	\date{}
	\maketitle
\begin{abstract}
		  A non-complete graph $G$ is said to be $t$-tough if for every vertex cut $S$ of $G$, the ratio of $|S|$ to the number of components of $G-S$ is at least $t$. The toughness $\tau(G)$ of the graph $G$ is the maximum value of $t$ such that $G$ is $t$-tough. A graph $G$ is said to be minimally $t$-tough if $\tau(G)=t$ and $\tau(G-e)<t$ for every $e\in E(G)$. 
          In 2003, Kriesell conjectured that every minimally $1$-tough graph contains a vertex of degree $2$. In 2018, Katona and Varga generalized this conjecture, asserting that every minimally $t$-tough graph contains a vertex of degree $\lceil 2t \rceil$.
		Recently, Zheng and Sun disproved the generalized Kriesell conjecture by constructing a family of $4$-regular graphs of even order. They also raised the question of whether there exist other minimally $t$-tough regular graphs that do not satisfy the generalized Kriesell conjecture. In this paper, we provide an affirmative answer by constructing a family of $4$-regular graphs of odd order, as well as a family of 6-regular graphs of order $3k+1~(k\geq 5)$.

			
\end{abstract}
	
	\section{Introduction}
	We consider finite simple graphs and use standard terminology and notation from \cite{Bondy} and \cite{West}.
	Let $G$ be a graph and let $c(G)$ denote the number of components
	of $G$. The graph $G$ is said to be $t$-tough if $t\cdot c(G-X)\le |X|$ for all $X\subseteq V(G)$ with $c(G-X)>1$, where $t$ is a nonnegative
	real number. The toughness $\tau(G)$ of the graph $G$ is the maximum value of $t$ such that $G$ is $t$-tough (taking $\tau(K_n)=\infty$ for
	all $n\ge 1$, where $K_n$ denotes the complete graph of order $n$). 
	The concept of toughness of a graph was introduced by
	Chv\'{a}tal~\cite{Chvatal}.
	Clearly, every hamiltonian graph is $1$-tough, but the converse is not true. Chv\'{a}tal~\cite{Chvatal} proposed the following
	conjecture, which is known as Chv\'{a}tal’s toughness conjecture. 
	\begin{conjecture}\label{chvatal}
		There exists a constant $t_0$ such that every $t_0$-tough graph of order at least three is hamiltonian.
	\end{conjecture}
	
	Bauer, Broersma and Veldman~\cite{Bauer2} showed that $t_0\ge \frac{9}{4}$ if it exists. Conjecture~\ref{chvatal} has been confirmed for a number of
	special classes of graphs~\cite{Bauer1,Kabela,Ota,Xu}.  Though great efforts have been made, it remains open.
	
	In 1999, Broersma, Engsberg and Trommel~\cite{Broersma} investigated a critical case of toughness and introduced the definition of minimally $t$-tough graphs. A graph $G$ is said to be minimally $t$-tough if $\tau(G) =t$ and $\tau(G-e)<t$ for all $e\in E(G)$. Katona, Solt\'{e}sz and Varga~\cite{Katona1} proved in 2018 the following result, demonstrating that the class of minimally $t$-tough graphs is large.
	\begin{theorem}[Katona, Solt\'{e}sz and Varga~\cite{Katona1}]
		For every positive rational number $t$, any graph can be embedded as an induced subgraph into a minimally $t$-tough graph.
	\end{theorem} 
	
	The connectivity of $G$, written $\kappa(G)$, is the minimum size of a vertex set $S$ such that $G-S$ is disconnected or has only one vertex. 
    It is easy to see that every minimally $t$-tough, noncomplete graph is $\lceil2t\rceil$-connected from the definition of toughness.
Therefore the minimum degree of a $t$-tough noncomplete graph is at least $\lceil2t\rceil$.
	The following conjecture is inspired by a result of Mader~\cite{Mader}, which asserts that every minimally $k$-connected graph contains a vertex of degree exactly $k$.
    
	\begin{conjecture}[Kriesell \cite{Kaiser}]\label{kriesell}
		Every minimally $1$-tough graph has a vertex of degree $2$.
	\end{conjecture}

Some results on Conjecture~\ref{kriesell} can be found in~\cite{Katona0,Katona1,Varga,Ma2}. In 2018, Katona and Varga~\cite{Katona2} generalized Kriesell’s conjecture to minimally $t$-tough graphs for any positive real number $t$.

\begin{conjecture}[Generalized Kriesell’s conjecture, Katona and Varga~\cite{Katona2}]\label{Conj4}
	Every minimally $t$-tough graph has a vertex of degree $\lceil 2t\rceil$, where $t$ is a positive real number.
\end{conjecture}
 
 Katona and Varga~\cite{Katona2} confirmed Conjecture~\ref{Conj4} on split graphs. In addition, Ma, Hu and Yang~\cite{Ma} proved that Conjecture~\ref{Conj4} holds for minimally $3/2$-tough, claw-free graphs. Recently, Zheng and Sun~\cite{Zheng} disproved Conjecture~\ref{Conj4} by constructing a family of $4$-regular graphs of even order. In the same paper, they posed the following problem.
 \begin{problem}\label{Pb5}
 	Are there other minimally $t$-tough regular graphs that do not satisfy the generalized Kriesell conjecture?
 \end{problem}

In this paper, we provide an affirmative answer to Problem~\ref{Pb5} by constructing an infinite family of $4$-regular graphs of odd order, which is a class of counterexamples for Conjecture~\ref{Conj4} distinct from the construction by Zheng and Sun~\cite{Zheng}. Moreover, for any integer $k\ge 5$, we construct a $6$-regular counterexample for Conjecture~\ref{Conj4}, which is a minimally $\frac{2k}{k-1}$-tough graph of order $3k+1.$ 

\begin{theorem}\label{main}
For every integer $k\ge 3$, there exists a $4$-regular minimally $\frac{k+1}{k-1}$-tough graph of order $2k+1$.
\end{theorem}

Clearly, for $k\ge 5$, the graph construted in Theorem~\ref{main} does not satisfy Conjecture~\ref{Conj4}, because $\lceil\frac{2(k+1)}{k-1}\rceil=3<4$.

\begin{figure}[!ht]
\centering
  \begin{tikzpicture}[scale = 0.6,main_node/.style={circle,fill=black,minimum size=0.68em,inner sep=0pt}]
			
			\foreach \i in {0,1,2,3,4,5,6,7,8,9,10} {
				\node[main_node] (\i) at (360/11*\i:4) {};  
			}
			
			\path[draw, black]
			(0) edge (1) 
			(1) edge (2) 
			(2) edge (3) 
			(3) edge (4) 
			(4) edge (5) 
			(5) edge (6) 
			(6) edge (7) 
			(7) edge (8) 
			(8) edge (9) 
			(9) edge (10) 
			(10) edge (0) 
			(0) edge (3) 
			(1) edge (4) 
			(2) edge (5) 
			(3) edge (6) 
			(4) edge (7) 
			(5) edge (8) 
			(6) edge (9) 
			(7) edge (10) 
			(8) edge (0) 
			(9) edge (1) 
			(10) edge (2);
			\node[below] at (0,-4.5){(a)};
		\end{tikzpicture}
 \hspace{1.5cm}
  \begin{tikzpicture}[scale = 0.6,main_node/.style={circle,fill=black,minimum size=0.68em,inner sep=0pt}]
%
\foreach \i in {0,1,2,3,4,5,6,7,8,9} {
				\node[main_node] (\i) at (360/10*\i:4) {};  
			}

			
			\path[draw, black]
(0) edge node {} (9) 
(9) edge node {} (8) 
(8) edge node {} (7) 
(7) edge node {} (6) 
(6) edge node {} (5) 
(5) edge node {} (4) 
(4) edge node {} (3) 
(3) edge node {} (2) 
(2) edge node {} (1) 
(0) edge node {} (1) 
(0) edge node {} (4) 
(1) edge node {} (5) 
(2) edge node {} (6) 
(3) edge node {} (7) 
(4) edge node {} (8) 
(5) edge node {} (9) 
(6) edge node {} (0) 
(7) edge node {} (1) 
(8) edge node {} (2) 
(9) edge node {} (3) 
(0) edge node {} (2) 
(1) edge node {} (3) 
(2) edge node {} (4) 
(3) edge node {} (5) 
(4) edge node {} (6) 
(5) edge node {} (7) 
(6) edge node {} (8) 
(7) edge node {} (9) 
(8) edge node {} (0) 
(9) edge node {} (1) ;

\node[below] at (0,-4.5){(b)};
		\end{tikzpicture}

\caption{The graphs constructed in the proofs of Theorems~\ref{main} and \ref{main2}}\label{Fig1}
\end{figure}
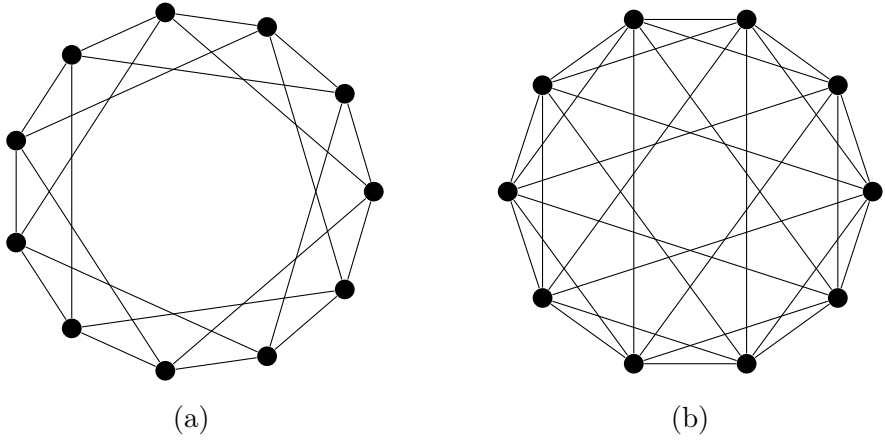

Due to Theorem~\ref{main}, as well as the counterexamples constructed by Zheng and Sun~\cite{Zheng}, all of which are $4$-regular graphs, it is natural to ask whether there exist counterexamples that are not $4$-regular. We provide a positive answer to this question in the following.

\begin{theorem}\label{main2}
    For every integer $k\ge 3$, there exists a $6$-regular minimally $\frac{2k}{k-1}$-tough graph of order $3k+1$.
\end{theorem}

It is easy to check that the graph construted in Theorem~\ref{main2} does not satisfy Conjecture~\ref{Conj4} for $k\ge 5$, since $\lceil\frac{4k}{k-1}\rceil=5<6.$


	\section{Proofs of Theorems~\ref{main} and \ref{main2}}
	

\begin{proof}[ {\bf Proof of Theorem~\ref{main}}]
    
Let $C$ be a cycle with vertices $v_1,v_2,\ldots,v_{2k+1}$ in order.
Let $G$ be the graph obtained from $C$ by adding the edges $v_iv_{i+3}$ with $1\le i\le 2k+1$,  where every subscript is understood to be modulo $2k+1$; see Figure~\ref{Fig1}(a). We will show that $G$ is minimally $\frac{k+1}{k-1}$-tough.
    
	We first prove that $\tau(G)= \frac{k+1}{k-1}$. Let $S$ be an arbitrary vertex cut of $G,$ we will see $\frac{|S|}{c(G-S)}\ge \frac{k+1}{k-1}$. 
    
	Suppose that $|S|\ge k+1.$ 
    Since $G$ is hamiltonian, $\alpha(G)\le \lfloor\frac{2k+1}{2}\rfloor=k$. If $\alpha(G)=k$, by symmetry we may assume that $\{v_1,v_3,\dots,v_{2k-1}\}$ is a maximum independent set of $G$, which contradicts $v_1v_{2k-1}\in E(G)$. So we have $\alpha(G)\le k-1$, and hence $c(G-S)\le k-1$. Therefore, $\frac{|S|}{c(G-S)}\ge \frac{k+1}{k-1}$, as desired.
    
	Now we assume that $|S|\le k$. If there exist two integers $i,j$ such that $\{v_i,v_{i+1},v_j,v_{j+1}\}\subseteq S$, as $G$ has a Hamilton cycle $C$, we have $c(G-S)\le |S|-2$ and hence
	\[
	\frac{|S|}{c(G-S)}\ge \frac{|S|}{|S|-2}\ge \frac{k+1}{k-1}.
	\]
	If there exists exactly one integer $i$ such that $\{v_i,v_{i+1}\}\subseteq S$, as $v_{i-1}v_{i+2}\in E(G)$, we have $c(G-S)\le s-2$ and so
	\[
	\frac{|S|}{c(G-S)}\ge \frac{|S|}{|S|-2}\ge \frac{k+1}{k-1}.
	\]
	Next, we assume that there exists no integer $i$ such that $\{v_i,v_{i+1}\}\subseteq S$, that is, $S$ is an independent set of $C$. Since $|C|=2k+1$ and $|S|\le k$, $C-S$ has a component, say $A$, of order at least two. Let $V(A)=\{v_t,v_{t+1},\ldots,v_{t+r}\}.$ Since $\{v_{t-2}v_{t+1},v_{t+r-1}v_{t+r+2}\} \subseteq E(G)$, we have $c(G-S)\le |S|-2$ and hence
	\[
	\frac{|S|}{c(G-S)}\ge \frac{|S|}{|S|-2}\ge \frac{k+1}{k-1}.
	\]
	Therefore, $\tau(G)\ge \frac{k+1}{k-1}$, as desired. 
    
	Now setting $X=\{v_{2i+1}\,|\,0\le i\le k\},$ one gets
	$\tau(G)\le \frac{|X|}{c(G-X)}=\frac{k+1}{k-1},$ 
    and hence $\tau(G)=\frac{k+1}{k-1}.$
	
	Secondly, we prove that $\tau(G-e)<\frac{k+1}{k-1}$ for any $e\in E(G)$. By symmetry, we only need to consider the following two cases.
\begin{wst}
\item[$\bullet$]
    $e=v_1v_2$. Let $S=\{v_{2i+1}\,|\,1\le i\le k\}$, we have $|S|=k$ and $c((G-e)-S)=k-1$. Thus, $\tau(G-e)\le \frac{k}{k-1}<\frac{k+1}{k-1}$.
  
\item[$\bullet$] 
    $e=v_1v_4$. Consider $U=\{v_2\}\cup\{v_{2i+1}\,|\,2\le i\le k\}$, we get $|U|=k$ and $c((G-e)-U)=k-1$. Thus, $\tau(G-e)\le \frac{k}{k-1}<\frac{k+1}{k-1}$.
\end{wst}
This completes the proof of Theorem~\ref{main}.
\end{proof}

  \begin{proof}[\bf Proof of Theorem~\ref{main2}]
    Let $C$ be a cycle with vertices $v_1,v_2,\ldots,v_{3k+1}$ in order. Let $G$ be the graph obtained from $C$ by adding the edges $v_iv_{i+2},v_iv_{i+4}$ with $1\le i\le 3k+1$, where every subscript is understood to be modulo $3k+1$; see Figure~\ref{Fig1}(b). 
    We will show that $G$ is minimally $\frac{2k}{k-1}$-tough.

    We first prove that $\tau(G)=\frac{2k}{k-1}$, and the following claim is needed.
    \begin{claim}\label{C1}
        $\alpha(G)\le k-1$.
    \end{claim}
    To the contrary, suppose that $I$ is an independent set of $G$ with $|I|=k.$ Without loss of generality, we may assume that $I=\{v_{i_1}, v_{i_2},\ldots, v_{i_k}\}$ with $1=i_1<i_2<\cdots <i_k$.  Since  $C$ is a Hamilton cycle of $G,$ $C-I$ has $k$ components. Note that $v_i$ is adjacent to each of $\{v_{i+2},v_{i+4}\}$  for $1\leq i\leq 3k+1$. This implies that each component of $C-I$  has order either at least four or two. If $C-I$ has a component of order at least four, then $|G|\ge 3k+2$, a contradiction. Hence, each component of $C-I$ has order two. However, it follows that $|G|=3k$, a contradiction. 
    This proves Claim~\ref{C1}.

     Let $S$ be an arbitrary vertex cut of $G,$ we aim to show that $\frac{|S|}{c(G-S)}\geq \frac{2k}{k-1}$. One gets $c(G-S)\leq k-1$ from Claim~\ref{C1}. 
    If $|S|\geq 2k$, then 
    \begin{align*}
        \frac{|S|}{c(G-S)}\geq \frac{2k}{c(G-S)}\geq \frac{2k}{k-1}.
    \end{align*}   
Now we assume that $|S|\leq 2k-1$. Let $H_1,H_2,\ldots, H_t$ be the components of $C[S],$ we proceed by considering the subsequent two cases. 

{\bf Case 1.} $\min\{|H_i| : 1\le i\le t\}=1.$ 

Let $\ell$ be the number of trivial components of $C[S]$. Without loss of generality, we assume that $H_1=\{v_{i_1}\}$, $H_2=\{v_{i_2}\},\ldots, H_{\ell}=\{v_{i_{\ell}}\}$, where $1=i_1<i_2<\cdots <i_{\ell}$.  For convenience, we denote $S'=\{v_{i_1},v_{i_2},\ldots,v_{i_{\ell}}\} $ and $S^*=S\setminus
  S'$.  Then 
  $$
  C-\{v_{i_r}\,|\,1\le r\le \ell\}+\{v_{i_{r-1}}v_{i_{r+1}}\,|\,1\le r\le \ell\}
$$ 
is a Hamilton cycle of $G-S'$. Clearly, 
 \begin{flalign*}
     c(G-S)\le c(C-S'-S^*)\leq \frac{|S^*|}{2} =\frac{|S|-\ell}{2}.
 \end{flalign*}
 If $l\geq 2$, then combining the fact that $|S|\le 2k-1$, we have
 \begin{flalign*}
     \frac{|S|}{c(G-S)}\geq \frac{2|S|}{|S|-\ell}\geq \frac{2|S|}{|S|-2}\geq \frac{2k}{k-1}.
 \end{flalign*}
 
Thus we assume that $\ell=1$. 
     If there exists a component of $C[S]$ with order at least three, then $c(G-S) \leq c(C-S\setminus \{v_1\})\leq \frac{|S|-2}{2}$. This implies that 
      \begin{flalign*}
       \frac{|S|}{c(G-S)}\geq \frac{2|S|}{|S|-2}\geq \frac{2k}{k-1}.
    \end{flalign*}
    We now assume that $|H_i| =2$ for each $2\le i\le t.$ Note that $C'=v_{3k+1}v_2\overrightarrow{C}v_{3k+1}$ is a Hamilton cycle of $G-v_1$, where  $v_2\overrightarrow{C}v_{3k+1}$ is the segment $v_2v_{3}\ldots v_{3k+1}$ in $C$. This implies that 
    $$c(C'-S)\leq \frac{|S\setminus\{v_1\}|}{2}=\frac{|S|-1}{2}.$$ 
    If each component of $C-S$ has order one, 
    as $|S|=2t-1\le 2k-1$, we have $|G|=3t-1<3k+1$, a contradiction. It follows that $C-S$ has a component, say $A$, of order at least two.
    Let $V(A)=\{v_i,v_{i+1},\dots,v_{i+j}\}$, 
    then either $v_{i-1}\neq v_1$ or $v_{i+j+1}\neq v_1$.
    Since $\{v_{i+j-1}v_{i+j+3},v_{i-3}v_{i+1}\}\subseteq E(G)$, we have 
    \[
    c(G-S)\le c(C-S\setminus \{v_1\})-1=\frac{|S|-3}{2}.
    \]
    Thus, combining the fact that $|S|\le 2k-1$, we get
    \begin{flalign*}
       \frac{|S|}{c(G-S)}\geq \frac{2|S|}{|S|-3}\geq \frac{2k}{k-1}.
    \end{flalign*}

{\bf Case 2.} $\min\{|H_i| : 1\le i\le t\}\ge 2.$ 

If there exists a component of $C[S]$ with order three, then $c(G-S)\leq \frac{|S|-3}{2}$ and hence
    \begin{flalign*}
       \frac{|S|}{c(G-S)}\geq \frac{2|S|}{|S|-3}\geq \frac{2k}{k-1}.
    \end{flalign*}
Suppose that $\max\{|H_i| : 1\le i\le t\}\ge 4.$ Then $c(G-S)\leq \frac{|S|-2}{2}$ and hence
    \begin{flalign*}
       \frac{|S|}{c(G-S)}\geq \frac{2|S|}{|S|-2}\geq \frac{2k}{k-1}.
    \end{flalign*}
We then assume that $|H_i| =2$ for each $1\le i\le t.$ Similar as the above argument, we may assert that there exists a component of $C-S$  of order at least two. Otherwise, as $|S|\le 2k-1$, we have $|G|=\frac{3|S|}{2}<3k+1$, a contradiction. It follows that 
$c(G-S)\leq c(C-S)-1=\frac{|S|-2}{2}.$
  Thus,  
  \begin{flalign*}
       \frac{|S|}{c(G-S)}\geq \frac{2|S|}{|S|-2}\geq \frac{2k}{k-1}.
    \end{flalign*}
 Therefore $\tau(G)\geq \frac{2k}{k-1}$, as desired.
 
  Let $X=\{v_1,v_3\}\cup \left(\bigcup\limits_{t=2}^k\{v_{3t-1},v_{3t}\}\right)$. This gives 
$$\tau(G)\leq \frac{|X|}{c(G-X)}=\frac{2k}{k-1},$$
and hence, we conclude that $\tau(G)=\frac{2k}{k-1}.$

  We now prove that $\tau(G-e)<\frac{2k}{k-1}$ for each $e\in E(G)$. By symmetry, it suffices to consider the following three cases.%
\begin{wst}
\item[$\bullet$]
    $e=v_1v_2$. Let $S=\{v_3,v_{3k-1},v_{3k+1}\}\cup \left(\bigcup\limits_{t=2}^{k-1}\{v_{3t-1},v_{3t}\}\right)$, we have $|S|=2k-1$ and $c((G-e)-S)=k-1$. Thus, $\tau(G-e)\le \frac{2k-1}{k-1}<\frac{2k}{k-1}$.
  
\item[$\bullet$] 
    $e=v_1v_3$. Consider $S=\{v_2,v_{3k-2},v_{3k},v_{3k+1}\}\cup \left(\bigcup\limits_{t=2}^{k-1}\{v_{3t-2},v_{3t-1}\}\right) $, we get $|S|=2k$ and $c((G-e)-S)=k$. Hence $\tau(G-e)\le \frac{2k}{k}<\frac{2k}{k-1}$.
   
\item[$\bullet$]
    $e=v_1v_5$. Setting $S=\{v_2,v_{3},v_{3k+1}\}\cup \left(\bigcup\limits_{t=2}^{k}\{v_{3t-2},v_{3t}\}\right) $ gives  $|S|=2k+1$ and $c((G-e)-S)=k$. Thus, $\tau(G-e)\le \frac{2k+1}{k}<\frac{2k}{k-1}$.
 \end{wst}   
    Therefore, $G$ is a minimally $\frac{2k}{k-1}$-tough graph.
	This completes the proof of Theorem~\ref{main2}.
  \end{proof} 
	
\section*{Acknowledgement}  The authors are grateful to Professor Xingzhi Zhan for his constant support and guidance. This research  was supported by the NSFC grant 12271170 and Science and Technology Commission of Shanghai Municipality (STCSM) grant 22DZ2229014.

\section*{Declaration}

	\noindent$\textbf{Conflict~of~interest}$
	The authors declare that they have no known competing financial interests or personal relationships that could have appeared to influence the work reported in this paper.
	
	\noindent$\textbf{Data~availability}$
	Data sharing not applicable to this paper as no datasets were generated or analysed during the current study.

\end{document}